\documentclass{article}

\usepackage{arxiv}

\usepackage[utf8]{inputenc} 
\usepackage[T1]{fontenc}    
\usepackage{hyperref}       
\usepackage{url}            
\usepackage{booktabs}       
\usepackage{amsfonts}       
\usepackage{nicefrac}       
\usepackage{microtype}      
\usepackage{lipsum}
\usepackage{algorithm}
\usepackage{algorithmic}
\usepackage{amsmath}
\usepackage{amssymb}
\usepackage{epstopdf}
\usepackage{enumitem}
\usepackage{float}
\usepackage{graphicx}
\usepackage{lineno,hyperref}
\usepackage{lscape}
\usepackage{morefloats}
\usepackage{multicol}
\usepackage{multirow}
\usepackage{subfigure}
\usepackage{ulem}
\usepackage{longtable}
\usepackage{multirow}
\usepackage{array}
\usepackage[table]{xcolor}
\usepackage{threeparttable}

\newcommand{\change}[1]{\textcolor{black}{#1}}

\title{On the use of dynamic mode decomposition for time-series forecasting of ships operating in waves 
}

\author{
  Andrea Serani$^{1,*}$, Paolo Dragone$^1$, Frederick Stern$^2$, Matteo Diez$^1$\\
  $^1$CNR-INM, National Research Council, Institute of Marine Engineering, Rome, Italy  \\
  $^2$IIHR--The University of Iowa, Iowa City, USA\\
  $^*$corresponding author: \texttt{andrea.serani@cnr.it} \\
}

\begin{document}
\maketitle

\begin{abstract}
In order to guarantee the safety of payload, crew, and structures, ships must exhibit good seakeeping, maneuverability, and structural-response performance, also when they operate in adverse weather conditions.
In this context, the availability of forecasting methods to be included within model-predictive control approaches may represent a decisive factor.
Here, a data-driven and equation-free modeling approach for forecasting of trajectories, motions, and forces of ships in waves is presented, based on dynamic mode decomposition (DMD). DMD is a data-driven modeling method, which provides a linear finite-dimensional representation of a possibly nonlinear system dynamics by means of a set of modes with associated frequencies. Its use for ship operating in waves has been little discussed and a systematic analysis of its forecasting capabilities is still needed in this context.
Here, a statistical analysis of DMD forecasting capabilities is presented for ships in waves, including standard and augmented DMD. The statistical assessment uses multiple time series, studying the effects of the number of input/output waves, time steps, time derivatives, along with the use of time-shifted copies of time series by the Hankel matrix. The assessment of the forecasting capabilities is based on four metrics: normalized root mean square error, Pearson correlation coefficient, average angle measure, and normalized average minimum/maximum absolute error. Two test cases are used for the assessment: the course keeping of a self-propelled 5415M in irregular stern-quartering waves and the turning-circle of a free-running self-propelled KRISO Container Ship in regular waves. Results are overall promising and show how state augmentation (using from four to eight input waves, up to two time derivatives, and four time-shifted copies) improves the DMD forecasting capabilities up to two wave encounter periods in the future. Furthermore, DMD provides a method to identify the most important modes, shedding some light onto the physics of the underlying system dynamics.
\end{abstract}


\section{Introduction}
%
%
Commercial and military vessels have to met International Maritime Organization (IMO) guidelines and North Atlantic Treaty Organization (NATO) Standardization Agreements (STANAG) in order to guarantee the safety of structures, payload, and crew in adverse weather conditions, providing good performance for seakeeping, maneuverability, and structures. In this regard, the prediction capability of ship performance in real seas conditions, along with the understanding of the physics involved, is of upmost importance for marine operations, such as aircraft/helicopter landing, ship-to-ship cargo transfer, off-loading of small boats, and off-shore operations. Also due to climate changes (as their effects may include higher sea states and even more adverse weather, impacting the safety of payload, structure, crew), prediction capabilities of ships in real seas have attracted the attention of recent NATO Science and Technology Organization (STO) Applied Vehicle Technology (AVT) task groups, such as AVT-280 ``Evaluation of Prediction Methods for Ship Performance in Heavy Weather'' (2017-2019) and AVT-348 ``Assessment of Experiments and Prediction Methods for Naval Ships Maneuvering in Waves'' (2021-2023), where the focus was on the assessment of prediction methods for ship seakeeping and maneuvering in waves, respectively\change{, and AVT-351 ``Enhanced Computational Performance and Stability \& Control Prediction for NATO Military Vehicles'' (2021-2023) focusing on the stability and control of military vehicles operating in realistic conditions.} 

Recently, computational fluid dynamics (CFD) studies have demonstrated their maturity for the prediction of ship performance in waves, including the computational tools assessment and validation (versus experimental campaign) in extreme sea conditions \cite{van2020prediction,serani2021-OE}. Nevertheless, the computational cost associated with the CFD analysis is generally very high, especially if statistical convergence of relevant estimators is sought after and complex hydro-structural problems are investigated via high-fidelity fluid-structure interaction solvers \cite{diez2022experimental}. 

\begin{figure*}[t]
\centering
\includegraphics[width=0.8\textwidth]{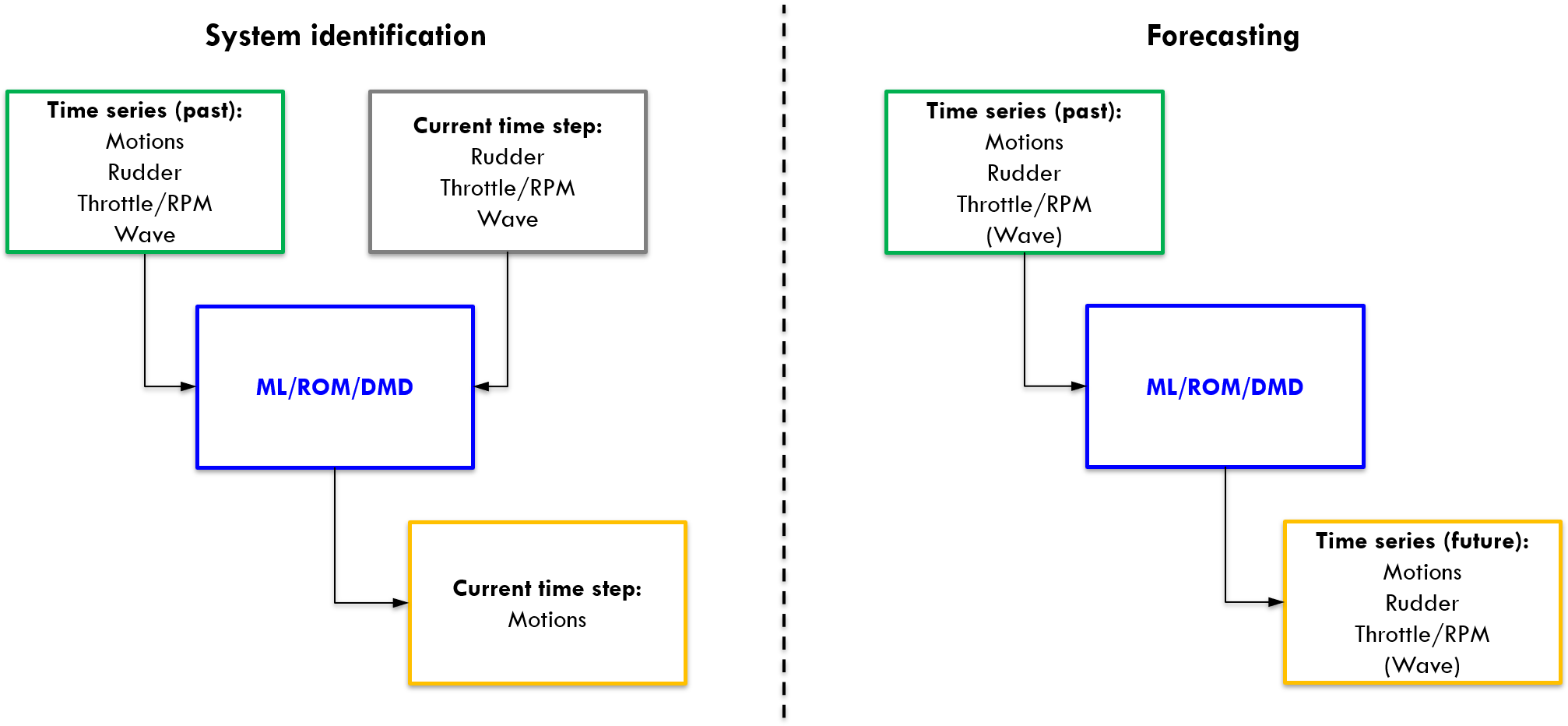}
\caption{System identification (left) and forecasting (right) conceptual schemes.}\label{fig:fcVSsi}
\end{figure*}

For these reasons, data-driven and equation-free approaches, such as machine learning (ML) methods and reduced-order models (ROM), may provide a viable option for the prediction of the ship performance in waves.
\change{
Furthermore and in parallel to prediction capabilities, the availability of accurate forecasting methods is essential to effective model-predictive control approaches, which may represent a decisive factor for the safety of ships in extreme weather.
It should be noted that ML and ROM approaches may be applied both to (a) the prediction of the performance for new ship designs or untested conditions and (b) the forecasting of relevant variables (such as motions, trajectories, and loads) based on available data. The first task usually follows a system identification setting \cite{silva2022data} whereas the second task is based on forecasting methodologies \cite{diez2022-OEME}. The conceptual scheme for system identification and forecasting settings is shown in Fig. \ref{fig:fcVSsi}.
}

Several approaches have been proposed to predict/forecast ship motions, including classic techniques based on statistical principles, as also ML methods, e.g. linear regression approaches have good forecasting performance and low computational cost, which makes them suitable for real-time realizations \cite{AR}. ML methods are very promising due to their capability of performing nonlinear modeling without an a priori knowledge of the relationships between input and output variables, like real-life ship motions \cite{ANN,SVR,NN,PSO-LSTM,LSTM,LSTM-GPR}. In particular, recurrent-type neural networks (RNN), such as long-short term memory (LSTM) and gated recurrent units (GRU), have shown their good short-term prediction/nowcasting capability for ship maneuvering in waves \cite{dagostino2022-OEME}. 
Nevertheless, according to current discussions in the data science community, these methods require more hands-on experience to set them up and make them efficient. Furthermore, ML computational cost is generally higher than statistical methods. 

ROM are generally easier to interpret than ML approaches and could help shedding light on the physics involved. \cite{MCA} developed a ship-motion prediction algorithm based on the minor component analysis, providing better prediction accuracy than neural networks (NNs) and auto-regressive (AR) models. To overcome the effect of the non-linearity and non-stationarity of ship-motion data, some hybrid forecasting approaches have been developed based on the empirical mode decomposition or closely related approaches \cite{AR-EMD-SVR,EMD-ARG-LSTM,EMD-SVR,EMD-GPR,EMD-CNN-GRU,diez2022snh}, used to decompose the data into a couple of intrinsic mode functions (IMF) and a residue, where each IMF represents a specific frequency component contained in the original signal; then a prediction model (e.g., AR, support-vector regression, Gaussian process regression, or RNN) of each component is built, implementing short-term prediction and than aggregating the predictions of all decomposed components. 

Among others, the dynamic mode decomposition (DMD) is increasing its success, especially in the fluid dynamics community (e.g., \cite{rowley2009spectral,timur2020}), since it is capable of providing accurate assessments of the space-time coherent structures in complex flows and systems \cite{magionesi2018modal,pagliaroli2022proper}, shedding the light on the problem physics, also allowing short-term future estimates of the system's state \cite{andreuzzi2021dynamic}, which can be used for real-time prediction and control \cite{kutz2016dynamic}.
DMD is a dimensionality-reduction/reduced-order modeling method, which provides a set of modes with associated oscillation frequencies and decay/ growth rates \cite{schmid2010dynamic}. For linear systems, these modes/frequencies correspond to the linear normal modes/frequencies of the system. More generally, DMD modes/frequencies approximate eigenmodes and eigenvalues of the infinite-dimensional linear Koopman operator, providing a linear finite-dimensional representation of the (possibly nonlinear) system dynamics \cite{kutz2016multiresolution}. 

Forecasting capabilities of the DMD has been investigated in several fields: \cite{electric-load} studied a short-term load forecasting model that employs DMD in time-delay embedded coordinates, integrated with the classical statistical method of Gaussian process regression; \cite{traffic-flow} proposed a low-rank DMD and reported better results in short-term traffic flow prediction compared to AR moving average and LSTM; \cite{offshore-wind-farm} present a data-driven real-time system identification approach to forecasting, based on streaming DMD (sDMD), where the algorithm works with data streams by adjusting the DMD continuously as new data are available, including also delay states to the matrices of measurements. 

For the ship hydrodynamics community, DMD offers a complementary efficient method to equation-based system identification approaches, e.g. \cite{araki2012estimating,araki2019improved}. The efficiency of the method in this context stems from the finite dimensionality of the set of relevant state variables together with the simplicity of operations required to model the system dynamics (as opposed to more data/resource-consuming ML approaches). This offers opportunity for integration into digital twin platforms \cite{fonseca2022standards} for the data-driven modeling and prediction of ships in waves\change{, as well as its beneficial application to steady and transient control for ship operations (e.g., \cite{li2016receding,liang2021finite,he2021adaptive}).} 

A proof-of-concept on the use of the DMD for the analysis and forecast of the finite-dimensional set of trajectory/ motion/ force time histories of ships operating in waves was presented by the authors in \cite{diez2022-OEME}. \change{Similar approaches to applying the DMD to finite-dimensional variable sets may be found for} power grid load data \cite{dylewsky2020dynamic,mohan2018data}, financial trading strategies \cite{mann2016dynamic}, sales data \cite{vasconcelos2019dynamic}, and neural recordings \cite{brunton2016extracting}. 
\change{Nevertheless, the application of DMD to model and forecast ships operating in waves has been little discussed in the literature and a systematic analysis of its prediction and forecasting capabilities is still needed in this context.
}

\change{
The objective of the current study is to extend the authors previous work \cite{diez2022-OEME}, providing a systematic statistical assessment of DMD forecasting capabilities of the finite-dimensional set of trajectory/ motion/ force time histories of ships operating in waves. 
It may be noted that DMD is a general methodology, which can be used both to extract knowledge of a dynamical process via dimensionality reduction and modal analysis, and for data-driven modeling and forecasting based on available data. Forecasting capabilities are in general problem dependent, for they depend on the dynamic characteristics of the system, the number of  degrees of freedom, and the availability of data. Here, the focus is on the assessment of DMD forecasting capabilities for trajectories, motions, and forces of ships in waves. In order to provide a robust assessment of these capabilities, a statistical analysis is proposed and discussed for two test cases of ships in waves. System identification capabilities of DMD (and its methodological extensions) is beyond the scope of the current work and will be addressed in future studies.
}

\change{Here}, DMD is augmented including time-series derivatives and time-shifted copies in the system's state matrix and a parametric analysis is conducted, varying the number of input and output waves, time derivatives (up to the fourth) and time-shifted copies (up to half wave). Numerical results are assessed based on four metrics, namely the normalized root mean square error, the Pearson correlation coefficient, the average angle measure, and the normalized average minimum/maximum absolute error. Results are presented and discussed for the course keeping of a free-running naval destroyer (5415M) in irregular stern-quartering waves, using unsteady Reynolds-averaged Navier-Stokes (URANS) computations \cite{serani2021-OE} data and the free-running KRISO Container Ship (KCS), using experimental data from the University of Iowa IIHR wave basin \cite{stern2022snh}, focusing on starboard turning circle in regular waves. 

The remainder of the paper is organized as follows. Section \ref{sec:DMD} presents the DMD \change{and the augmented formulations used for the current study}. A brief description of the test cases is given in Section \ref{sec:cases}.
\change{Section \ref{sec:setup} describes the DMD setup and how the systematic study and statistical analyses for ships in waves are performed}. 
Section \ref{sec:results} presents and \change{discusses} the numerical results, and, finally, conclusions \change{on the use of DMD for ships in waves} and possible extension for future work are discussed in Section \ref{sec:conclusions}.

\section{Dynamic Mode Decomposition}\label{sec:DMD}
Consider a dynamical system described as \cite{kutz2016dynamic}
\begin{equation}\label{eq:sysdyn}
\dfrac{\mathrm{d}\mathbf{x}}{\mathrm{d}t}=\mathbf{f}(\mathbf{x},t;\mu),
\end{equation}
where $\mathbf{x}(t)\in\mathbb{R}^n$ represents the system's state at time $t$, $\mu$ contains the parameters of the system, and $\mathbf{f}(\cdot)$ represents its dynamics.
The state $\mathbf{x}$ is generally large, with $n\gg 1$ and can represents, for instance, the discretization of partial differential equations at a number of discrete spatial points, or multi-channel/multi-variable time series.

Considering $\mathbf{f}(\mathbf{x},t;\mu)$ as unknown, the DMD works with an equation-free perspective. Thus, only the system measurements are used to approximate the system dynamics and forecast the future states. Equation \ref{eq:sysdyn} is approximated by the DMD as a locally linear dynamical system defined as \cite{kutz2016dynamic}
\begin{equation}\label{eq:dmdsys}
\dfrac{\mathrm{d}\mathbf{x}}{\mathrm{d}t}=\mathcal{A}\mathbf{x}
\end{equation}
with solution
\begin{equation}\label{eq:dmdrec}
\mathbf{x}(t)=\sum_{k=1}^n \boldsymbol{\phi}_k \, q_k(t)=
\sum_{k=1}^n \boldsymbol{\phi}_k \, b_k \,\exp(\zeta_k t),
\end{equation}
where $\boldsymbol{\phi}_k$ and $\zeta_k$ are, respectively, the eigenvectors and the eigenvalues of the matrix $\mathcal{A}$, $q_k$ are the time-varying modal coordinates, and the coefficients $b_k$ are the modal coordinates of the initial condition $\mathbf{x}_0$ in the eigenvector basis, for which
\begin{equation}
\mathbf{b}=\boldsymbol{\Phi}^{-1}\mathbf{x}_0.
\end{equation}

Sampling the system every $\Delta t$, the time-discrete state can be expressed as $\mathbf{x}_k=\mathbf{x}(k\Delta t)$ with $k=1,...,m$, representing from now on known the system measurements. Consequently, the equivalent discrete-time representation of the system in Eq. \ref{eq:dmdsys} can be written as
\begin{equation}\label{eq:dmddsys}
\mathbf{x}_{k+1}=\mathbf{Ax}_k, \qquad \mathrm{with} 
\qquad \mathbf{A}=\exp(\mathcal{A}\Delta t)
\end{equation}

Arranging all the $m$ system measurements in the following two matrices
\begin{equation}\nonumber
\mathbf{X}=
\begin{bmatrix}
\mathbf{x}_k & \mathbf{x}_{k+1} & \dots & \mathbf{x}_{m-1}\\
\end{bmatrix},
\end{equation}
\begin{equation}\label{eq:XX'}
\mathbf{X}'=
\begin{bmatrix}
\mathbf{x}_{k+1} & \mathbf{x}_{k+2} & \dots & \mathbf{x}_{m}\\
\end{bmatrix},
\end{equation}
the matrix $\mathbf{A}$ in Eq. \ref{eq:dmddsys} can be constructed using the following approximation
\begin{equation}\label{eq:approxA}
\mathbf{A}\approx\mathbf{X}'\mathbf{X}^{\dag},
\end{equation}
where $\mathbf{X}^{\dag}$ is the Moore-Penrose pseudo-inverse of $\mathbf{X}$, which minimize $\|\mathbf{X}'-\mathbf{AX}\|_F$, where $\|\cdot\|_F$ is the Frobenius norm. 

The state-variable evolution in time can be approximated by the following modal expansion, as per Eq. \ref{eq:dmdrec}, 
\begin{equation}
\mathbf{x}(t)=\sum_{k=1}^n \boldsymbol{\varphi}_k \, q_k(t)=
\sum_{k=1}^n \boldsymbol{\varphi}_k \, b_k \,\exp(\omega_k t),    
\end{equation}
where $\boldsymbol\varphi_k$ are the eigenvectors of the approximated matrix $\bf A$ and $\omega_k=\text{ln}(\lambda_k)/\Delta t$, with $\lambda_k$ the eigenvalues of the same matrix \cite{kutz2016dynamic}.

From a general point of view, the DMD can be viewed as a method to compute the eigenvalues and eingevectors (modes) of a finite-dimensional linear model that approximates the infinite-dimensional linear Koopman operator \cite{kutz2016dynamic}, also known as the composition operator. 
Furthermore, due to the low dimensionality of data in the current context, Eq. \ref{eq:approxA} is computed directly, without the need of performing the singular value decomposition of $\bf X$ and projecting onto proper orthogonal decomposition modes \cite{kutz2016dynamic}. 


\subsection{Augmented DMD}
The DMD formulation can be extended through the so called state augmentation. In the following this is achieved by transforming the matrix $\bf A$ in the augmented matrix $\widehat{\bf A}$, where the notation ` $\widehat{\cdot}$ ' is used to indicate the augmented matrix. Specifically, three augmenting strategies are explored adding to the state vector: (1) time derivatives ${\mathrm{d}^i\mathbf{x}}/{\mathrm{d}t^i}$, (2) time-shifted copies (delayed states), and (3) both time derivatives and time-shifted copies.  

In the first case, defining the time-derivatives matrices as 
\begin{equation}\label{eq:DD'}
\mathbf{D}=
\begin{bmatrix}
\dfrac{\mathrm{d}\mathbf{X}}{\mathrm{d}t} \\ 
\\
\dfrac{\mathrm{d}^2\mathbf{X}}{\mathrm{d}t^2} \\ 
\\
\vdots \\
\\
\dfrac{\mathrm{d}^i\mathbf{X}}{\mathrm{d}t^i}\\
\end{bmatrix},
\qquad
\mathbf{D}'=
\begin{bmatrix}
\dfrac{\mathrm{d}\mathbf{X}'}{\mathrm{d}t} \\ 
\\
\dfrac{\mathrm{d}^2\mathbf{X}'}{\mathrm{d}t^2} \\ 
\\
\vdots \\
\\
\dfrac{\mathrm{d}^i\mathbf{X}'}{\mathrm{d}t^i}\\
\end{bmatrix},
\end{equation}
yields an augmented version of the matrices in Eq. \ref{eq:XX'} as
\begin{equation}\label{eq:dXX'}
\widehat{\mathbf{X}}=
\begin{bmatrix}
\mathbf{X} \\
\mathbf{D}\\
\end{bmatrix},
\qquad
\widehat{\mathbf{X}}'=
\begin{bmatrix}
\mathbf{X}' \\ 
\mathbf{D}'\\
\end{bmatrix}.
\end{equation}
\change{For the present study, time derivatives are evaluated via a backward finite difference scheme, so that only past states are needed.}

In the second case, defining the matrix of the time-shifted copies in the form of a Hankel matrix
\begin{equation}\nonumber
\mathbf{S}=
\begin{bmatrix}
\mathbf{x}_{k-1} & \mathbf{x}_{k} & \dots & \mathbf{x}_{m-2}\\
\mathbf{x}_{k-2} & \mathbf{x}_{k-1} & \dots & \mathbf{x}_{m-3}\\
\vdots & \vdots & \vdots & \vdots \\
\mathbf{x}_{k-s-1} & \mathbf{x}_{k-s} & \dots & \mathbf{x}_{m-s-1}\\
\end{bmatrix},
\end{equation}
\begin{equation}\label{eq:sts'}
\mathbf{S}'=
\begin{bmatrix}
\mathbf{x}_{k} & \mathbf{x}_{k+1} & \dots & \mathbf{x}_{m-1}\\
\mathbf{x}_{k-1} & \mathbf{x}_{k} & \dots & \mathbf{x}_{m-2}\\
\vdots & \vdots & \vdots & \vdots \\
\mathbf{x}_{k-s} & \mathbf{x}_{k-s+1} & \dots & \mathbf{x}_{m-s}\\
\end{bmatrix},
\end{equation}
yields an augmented version of the matrices in Eq. \ref{eq:XX'} as
\begin{equation}\label{eq:sXX'}
\widehat{\mathbf{X}}=
\begin{bmatrix}
\mathbf{X} \\
\mathbf{S}\\
\end{bmatrix},
\qquad
\widehat{\mathbf{X}}'=
\begin{bmatrix}
\mathbf{X}' \\ 
\mathbf{S}'\\
\end{bmatrix}.
\end{equation}

Finally, using both time derivatives and time-shifted copies
provides an augmented version of the matrices in Eq. \ref{eq:XX'} as
\begin{equation}\label{eq:stXX'}
\widehat{\mathbf{X}}=
\begin{bmatrix}
\mathbf{X} \\
\mathbf{D}\\
\mathbf{S}\\
\end{bmatrix},
\qquad
\widehat{\mathbf{X}}'=
\begin{bmatrix}
\mathbf{X}' \\ 
\mathbf{D}'\\
\mathbf{S}'\\
\end{bmatrix};
\end{equation}

\change{It is worth noting that, for the current study, the number of state and augmented variables is relatively low, as consequently also the dimensionality of the state matrix. Therefore, the complexity and computational cost associated with the overall modeling and forecasting process, including that of solving the matrix eigenproblem, is not a critical issue. Dimensionality/model-order reduction is not necessary, as in standard DMD. Conversely, the model order here is increased by the augmented states with the aim of gaining modeling and forecasting accuracy.}

\section{Test Cases Description}\label{sec:cases}
Two test cases are selected for the assessment of the DMD performance. Specifically, the CFD seakeeping simulations of a self-propelled destroyer-type vessel (5415M) in irregular stern-quartering waves and the EFD maneuvering data of a free-running container ship (KCS) turning circle in regular waves. Further details are given in the following subsections. 

\begin{figure}[t]
\centering
\includegraphics[width=0.5\columnwidth]{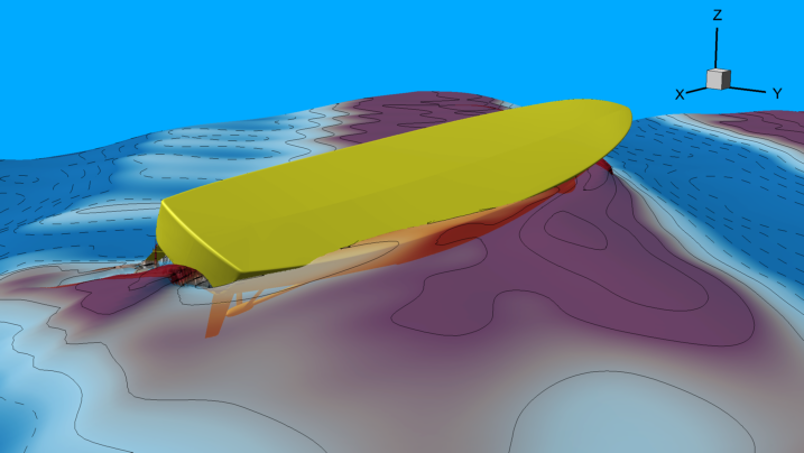}
\caption{CFD snapshot of the 5415M test case.}\label{fig:5415}
\end{figure}

\subsection{Course Keeping of the 5415M in Irregular Waves}
The hull form under investigation is the MARIN model 7967 which is equivalent to 5415M, used as test case for the NATO AVT-280 ``Evaluation of Prediction Methods for Ship Performance in Heavy Weather'' \cite{van2020prediction}. The model is self-propelled and kept on course by a proportional-derivative controller actuating the rudder angle. 

Course-keeping computations are based on the URANS code CFDShip-Iowa V4.5 \cite{huang2008-IJNMF}. CFD simulations are performed with propeller revolutions per minute fixed to the self-propulsion point of the model for the nominal speed, corresponding to $\mathrm{Fr} = 0.33$. The simulations are conducted in irregular long crested waves (following a JONSWAP spectrum), with nominal peak period $T_p = 9.2$\textit{s} and wave heading of 300\textit{deg} (see Fig. \ref{fig:5415}). The nominal significant wave height is equal to 7\textit{m}, corresponding to sea state 7 (high), according to the World Meteorological Organization definition. 
%
The six degrees of freedom rigid body equations of motion are solved to calculate linear and angular motions of the ship. A simplified body-force model is used for the propeller, which prescribes axisymmetric body force with axial and tangential components. The total number of grid points is about 45M. Further details can be found in \cite{serani2021-OE} and \cite{van2020prediction} where also potential flow computations and experimental data are presented and discussed. 

The state variables used for DMD are the ship six degrees of freedom (surge, sway, heave, roll, pitch, yaw), plus the rudder angle. Note that ship motions are in the carriage coordinate system, projected onto the ship axes.
\change{It may be noted that here the DMD is proposed for the solution of a data-driven modelling problem, using a reduced set of variables that can be reasonably measured onboard. The addition of wave field variables (using in practice, e.g., wave radar data) would be a natural extension of the current study and subject of future work.}

\subsection{Turning Circle of the KCS in Regular Waves}
\begin{figure}[t]
\centering
\includegraphics[width=0.5\columnwidth]{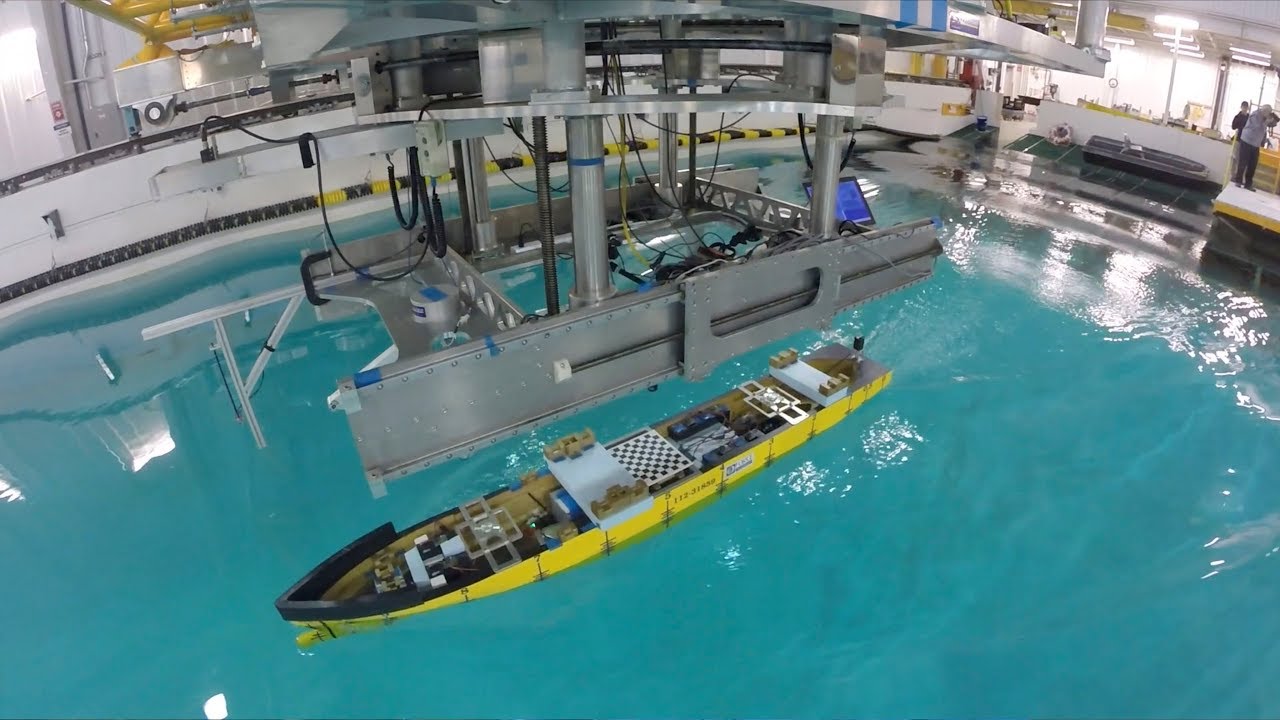}
\caption{Photograph of the IIHR wave basin with the KCS model.}\label{fig:KCS}
\end{figure}

The starboard turning circle of the free-running KCS in regular waves with constant rudder angle of 35\textit{deg} is taken as second test case. Data are taken from experiments conducted at the IIHR wave basin, which is shown in Fig. \ref{fig:KCS} and whose characteristics are given in \cite{sanada2021experimental}. The model length is $L=2.7$\textit{m} and the nominal speed corresponds to $\mathrm{Fr}=0.157$. The propeller RPM are fixed and provide the nominal speed in calm water. The regular wave parameters are $\lambda/L=1$ (wave length to ship-length ratio) and $H/\lambda=1/60$ (wave height to wave-length ratio).

The state variables used for DMD are the $x$, $y$, and $z$ coordinates (Earth-coordinate system) of a reference point placed amidships, \change{$u$, $v$, and $w$ components of the ship velocity (projected onto the ship axes), pitch and roll motions, drift angle, and propeller thrust and torque}. 
\change{As for the 5415M case, a reduced set of variables is used that can be reasonably measured onboard. The addition of wave field variables would be also in this case a natural extension of the current study and will be addressed in future work.}

\section{Numerical Setup and Metrics}\label{sec:setup}
Although the DMD is not a machine learning method in the strict sense, its data-driven nature allows for approaching DMD in a similar way to machine learning. 
Here the matrix $\bf A$ is constructed using observed (past) time histories, which are used as training set. The DMD is then used for the short-term prediction of system state time histories, which are compared against true observed (future) data used as test set. 

For each test case, different DMD setups are implemented by varying the dimension of the training set and the number of derivatives and/or time-shifted copies used to arrange the matrices of measurements. These solutions are investigated combining all possibilities (full factorial design of experiments). Dimension of training set are expressed in terms of encounter waves (or periods/oscillations of the phenomena) in order to have a faster/friendly physic interpretation. The test set is also expressed as numbers of future oscillations. Specifically, both 5415M and KCS time-histories are normalized with the time of one encounter waves, here equal to 32 time steps for both cases. Up to the fourth time derivative of the system states are used for the ADMD, where derivatives are calculated through a second-order accuracy backward finite difference scheme. For the time-shifted copies, up to 16 time shifts are used, that correspond to half encounter wave. The training sets are composed by a number of encounter waves from 1 up to 8 (in the past), whereas the test set goes till 4 encounter waves (in the future). 

Future state predictions/forecast are built using all modes/ frequencies. 
\change{
It is worth noting that in the results section, frequencies and eigenvectors are ordered by the
absolute value of the complex frequency imaginary part, in ascending order. Alternatively, frequencies and eigenvectors may be ordered by their modal participation or energy, in descending order, which is a natural choice if model-order reduction is sought after. This is not needed here as the number of state and augmented variables (and therefore the state matrix dimensionality) is relatively low. Therefore, the computational cost of the process is not a critical issue and a reduced-dimensionality/ order model is not necessary.
}

All variables are standardized, i.e., translated and scaled to have zero mean and unit variance. The prediction performance are investigated in a statistical way, considering 1001 random samples of the available data set, where one sample is composed by a number of input and a number of output wave (NIW and NOW) and the state variables for the DMD, along with derivatives and delay copies for the ADMD.
\begin{figure*}[!t]
\centering
\includegraphics[width=0.95\textwidth]{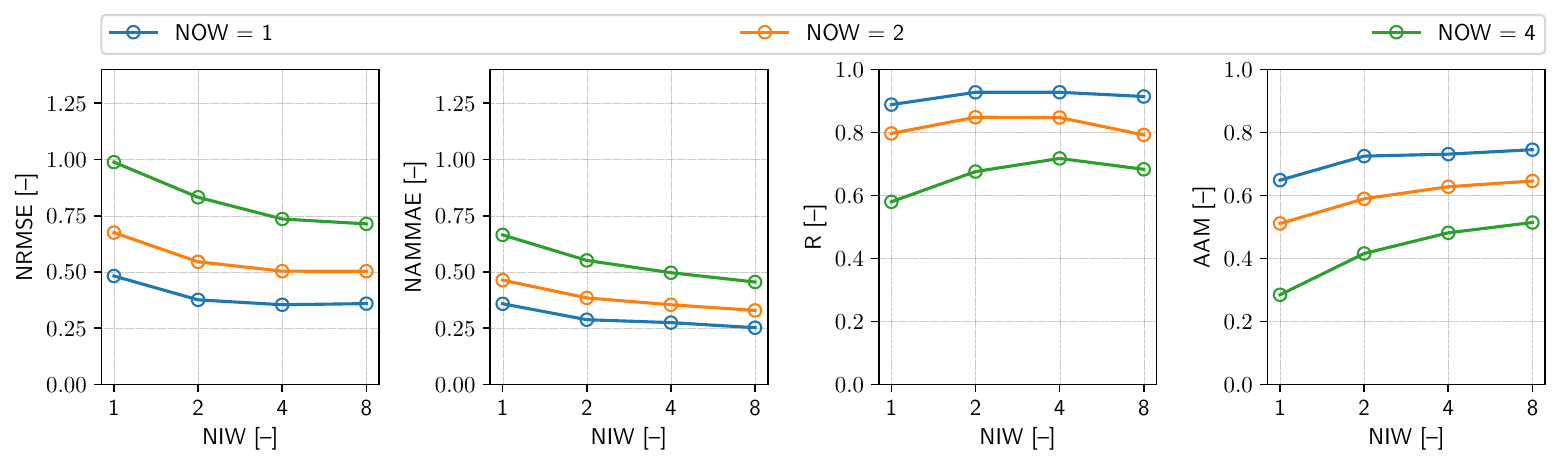}
\caption{DMD prediction performance on 5415M case; from left to right medians of NRMSE, NAMMAE, R, and AAM.}
\label{fig:DMDstd_5415M}
\end{figure*}
\begin{figure*}[!t]
\centering
\includegraphics[width=0.95\textwidth]{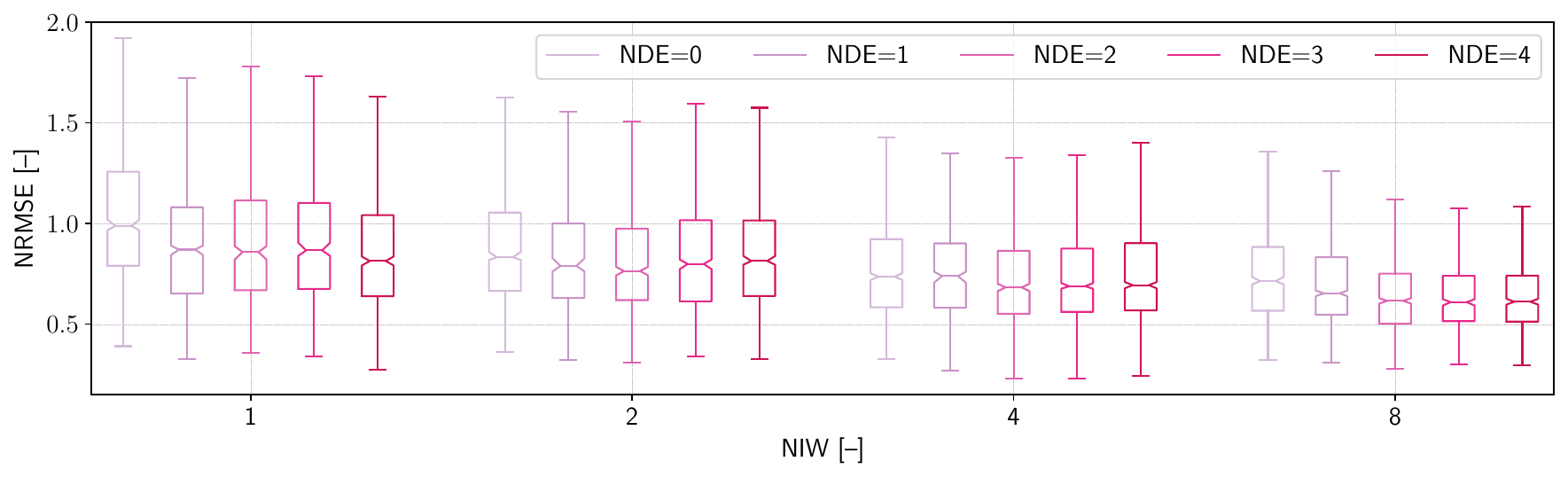}
\caption{Box plots of NRMSE values of ADMD with derivatives for 5415M test case.}\label{fig:DMDder_5415M}
\end{figure*}

Four evaluation metrics are used to assess the accuracy of the output sequences versus test time histories, namely the normalized root mean square error (NRMSE), the Pearson correlation coefficient (R), the average angle measure (AAM), and the normalized average minimum/maximum absolute error (NAMMAE). Metrics are averaged across the output variables to provide an overall assessment of the prediction accuracy.

The average NRMSE is defined as
\begin{equation}
\mathrm{NRMSE}=\frac{1}{n}\sum_{i=1}^n \sqrt{\frac{1}{m \sigma_{y_i}^2}\sum_{j=1}^m [x_i(t_j) - y_i(t_j)]^2}
\end{equation}  
where $x_i$ and $y_i$ represent respectively predicted and measured (test) values for the $i$-th variable; $\sigma_{y_i}$ represents the standard deviation of $y_i$. 

Similarly, the average Pearson correlation coefficient is evaluated as
\begin{equation}
R=\frac{1}{n(m-1)}  \sum_{i=1}^n
\sum_{j=1}^m \frac{[x_i(t_j)-\bar{x}_i]^2}{\sigma_{x_i}}
             \frac{[y_i(t_j)-\bar{y}_i]^2}{\sigma_{y_i}}
\end{equation}
where  $\bar{x}$ represents the time mean of $x$.

The AAM is a metric developed at NSWCCD in the 90s to quantify the accuracy of predicted time series when compared to actual measured time series \cite{hess2006feedforward}. Its averaged version reads
\begin{equation}
\mathrm{AAM}=\frac{1}{n} \sum_{i=1}^n \left\{
1-\frac{4}{\pi}
\left[ 
\frac{\sum_{j=1}^m d_i(t_j) |\alpha_i(t_j)|}{ \sum_{j=1}^m d_i(t_j)}
\right]\right\}
\end{equation}
where
\begin{equation}
\alpha_i(t_j) = \arccos \left[\frac{|x_i(t_j) + y_i(t_j)|}{\sqrt{2} d_i(t_j)} \right] 
\end{equation}
and
\begin{equation}
d_i(t_j) = \sqrt{x_i^2(t_j)+y_i(t_j)^2}  
\end{equation}

Finally, the NAMMAE metric is introduced \cite{diez2022snh} to provide a more engineering-oriented assessment of the prediction accuracy, based on minimum and maximum values of predicted and measured time series, as
\begin{eqnarray}\nonumber
\mathrm{NAMMAE} = \sum_{i=1}^n
\left\{  
         \left| \min_j{[x_i(t_j)]} - \min_j{[y_i(t_j)]} \right| \right. & \\
\left. + \left| \max_j{[x_i(t_j)]} - \max_j{[y_i(t_j)]} \right| \right\} \frac{1}{2 n \sigma_{y_i}}
\end{eqnarray}
\begin{figure*}[!t]
\centering
\includegraphics[width=0.95\textwidth]{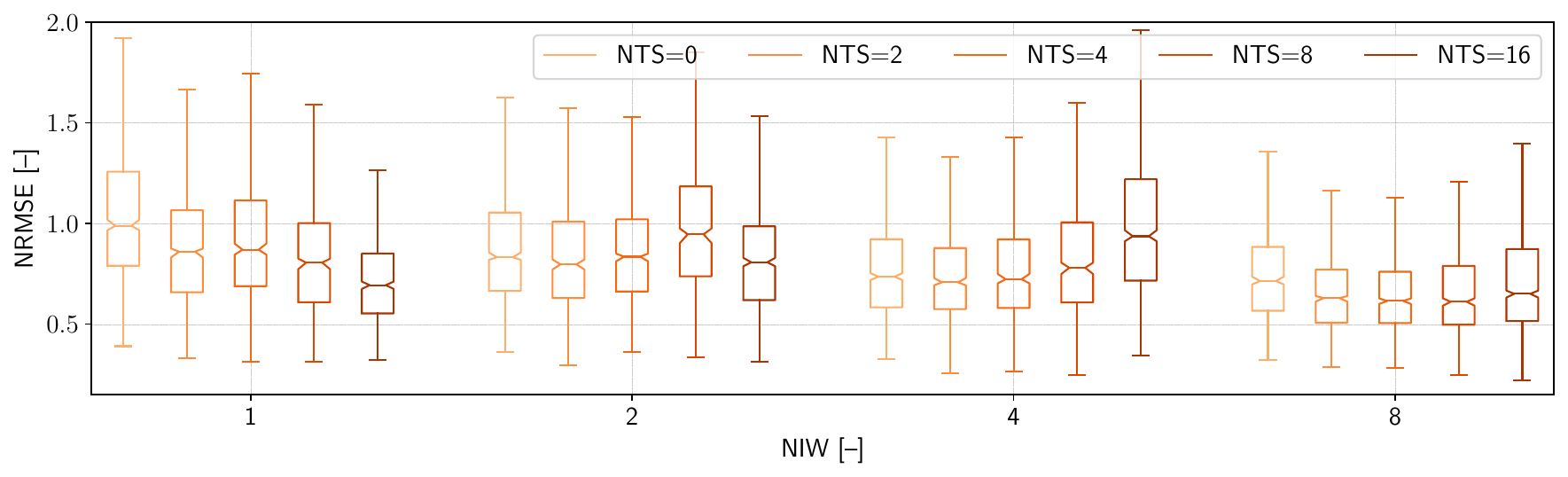}
\caption{Box plots of NRMSE values of ADMD with time-shifted copies for 5415M test case.}\label{fig:DMDshift_5415M}
\end{figure*}
\begin{figure*}[!t]
\centering
\includegraphics[width=1\textwidth]{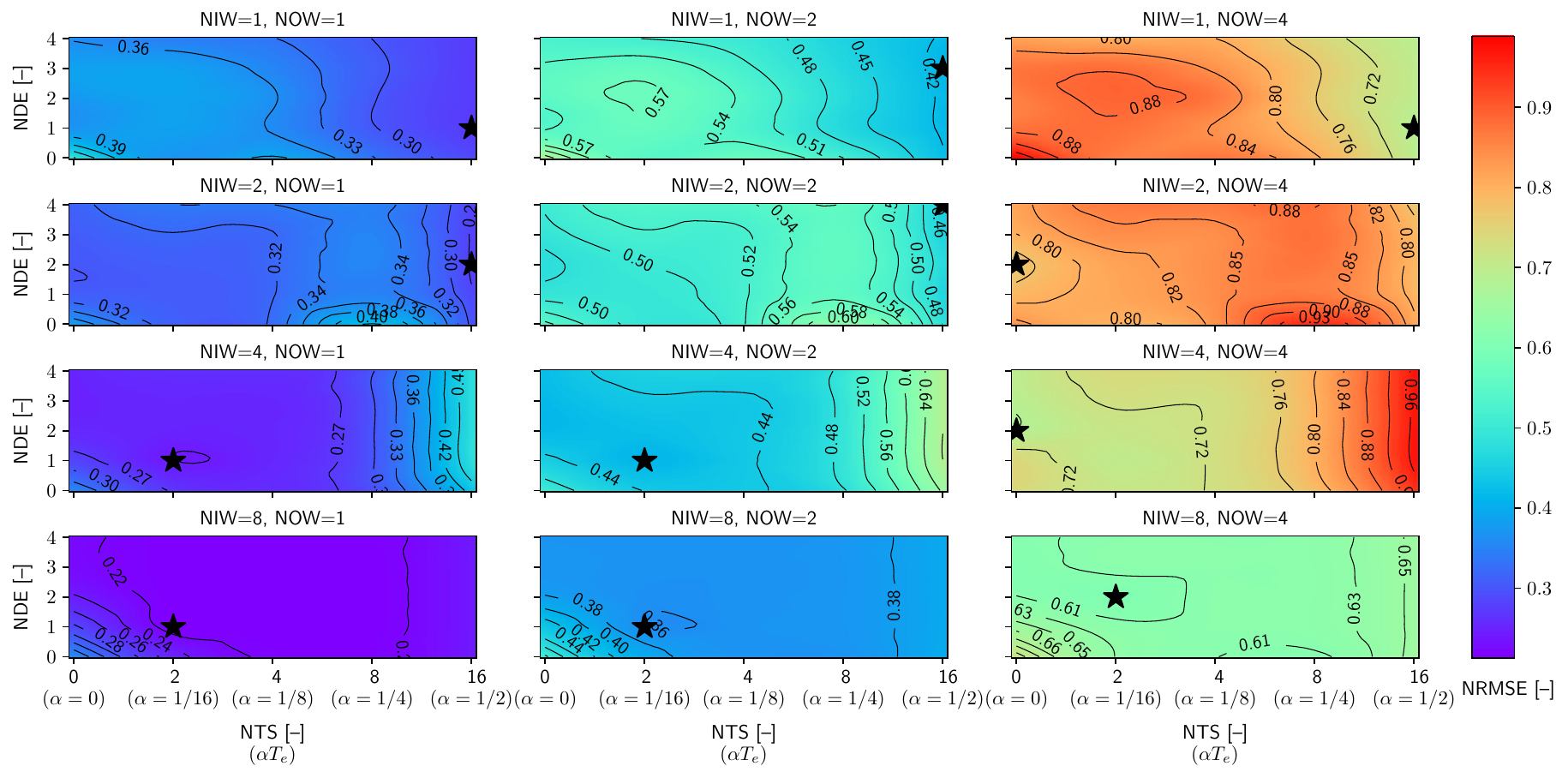}
\caption{Medians of NRMSE for the full factorial analysis of the DMD setups for 5415M test case.}\label{fig:DMDcplot_5415M}
\end{figure*}


\section{Results and Discussion}\label{sec:results}
The following subsections describe the results of the statistical analysis conditional to DMD setup for 5415M and KCS test cases, including examples of system states forecast. Finally, modal analysis discussion is also provided. It may be noted that, in the case that DMD provides eigenvalues with positive real part, that corresponds to a divergence condition (from a mathematical view point) of the system dynamics, a correction is imposed. Specifically, since it is known a priori that the system dynamics (from a physical view point, at least for the present test cases) is not diverging, if complex frequencies with positive real part are \change{found}, their real part is set to zero, in order to guarantee \change{stable} predictions of the future states.

\subsection{Forecasting of 5415M Course Keeping in Irregular Waves}

\begin{figure*}[!t]
\centering
\includegraphics[width=0.9\textwidth]{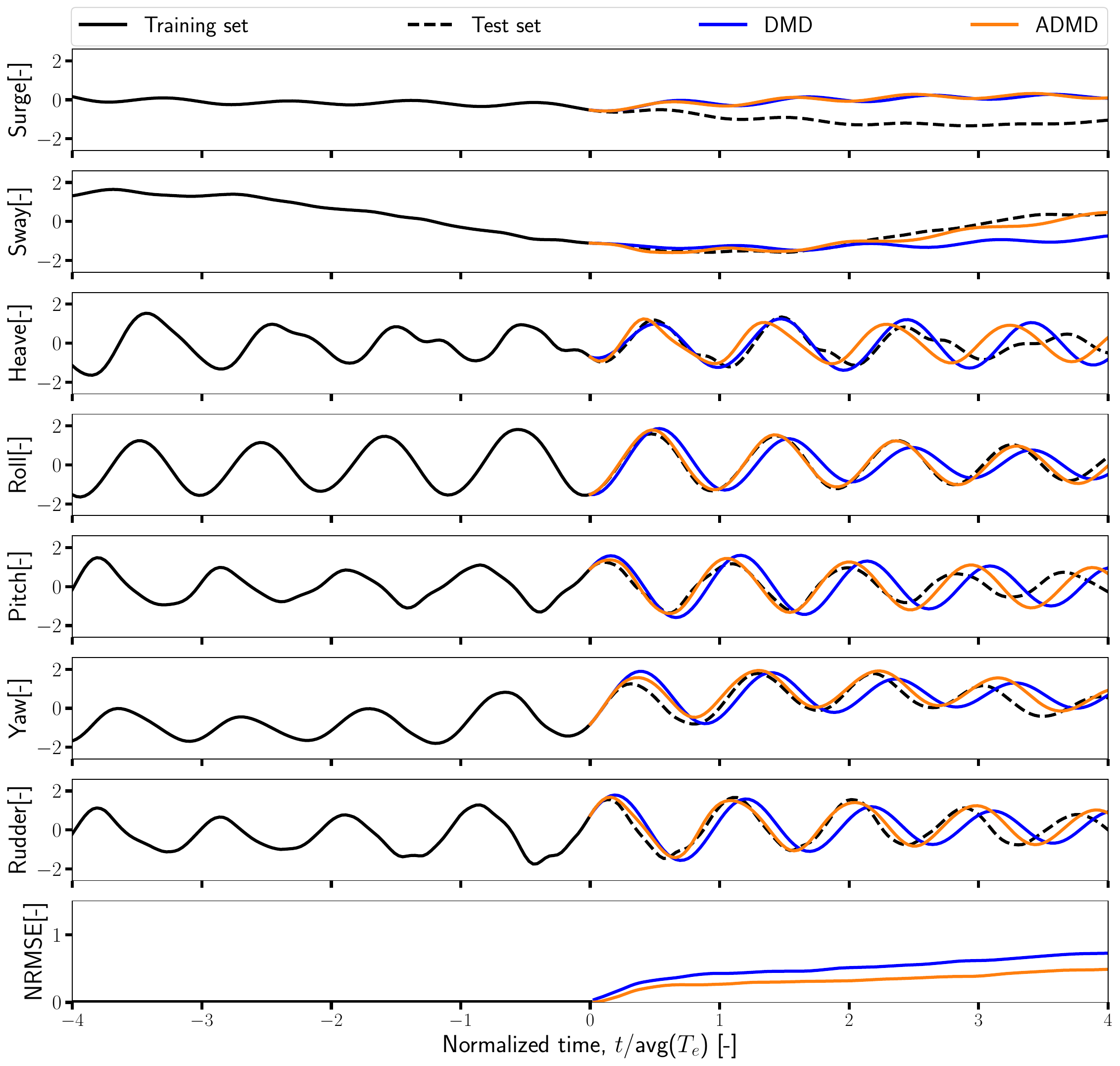}
\caption{Systems-state prediction comparison of DMD and best ADMD setup for a random time history of 5415M data set.}\label{fig:DMDvsADMD_5415M}
\end{figure*}
\begin{figure*}[!t]
\centering
\includegraphics[width=1\textwidth]{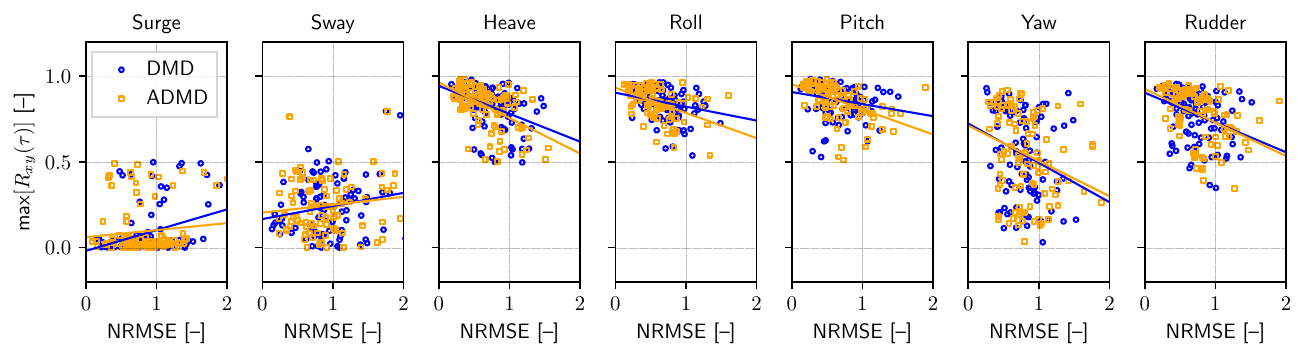}
\caption{Training/test sets maximum cross-correlation versus DMD and ADMD prediction error (NRMSE) for 5415M case.}\label{fig:maxCorr_5415M}
\end{figure*}
\begin{figure*}[!t]
\centering
\includegraphics[width=1\textwidth]{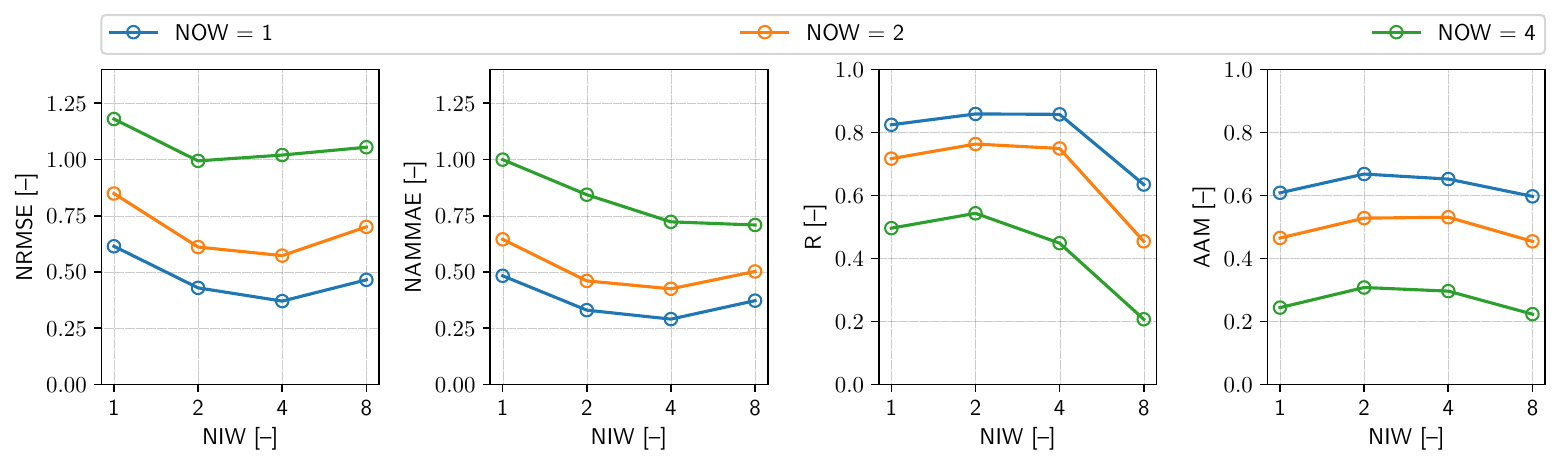}
\caption{DMD prediction performance on KCS case, in terms of medians of the evaluation metrics. From left to right: NRMSE, NAMMAE, R, and AAM.}\label{fig:DMDstd_KCS}
\end{figure*}
\begin{figure*}[!t]
\centering
\includegraphics[width=1\textwidth]{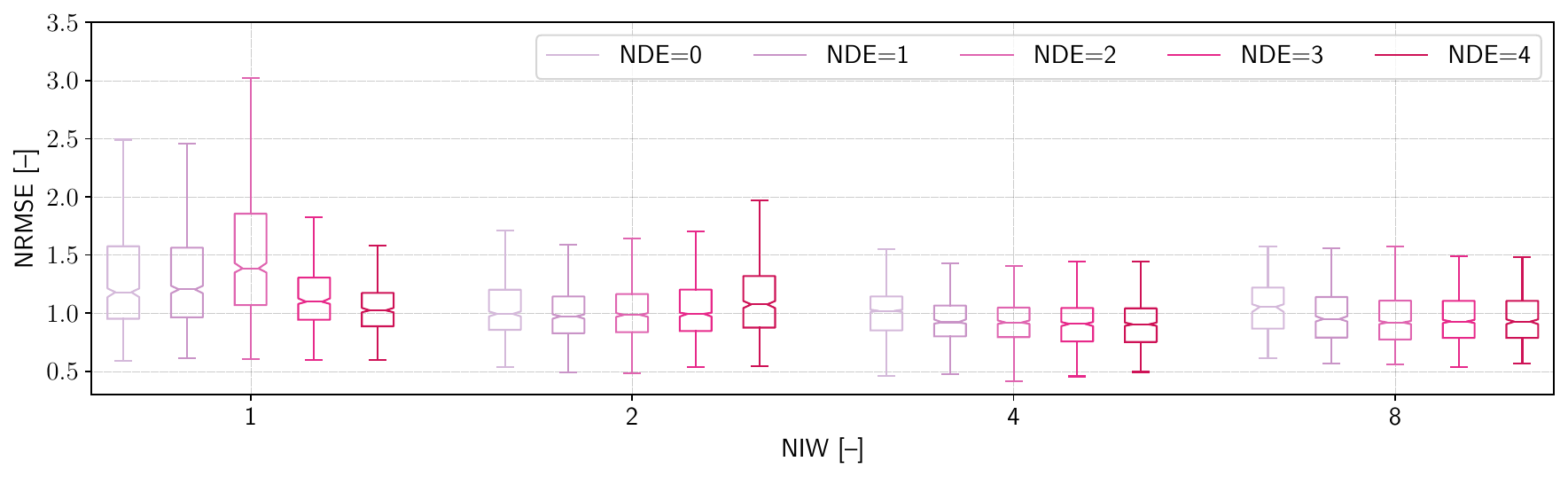}
\caption{Box plots of NRMSE values of ADMD with derivatives for KCS test case.}\label{fig:DMDder_KCS}
\end{figure*}

Figure \ref{fig:DMDstd_5415M} presents the prediction performance of the DMD standard formulation, showing the median values (of the statistic sample) for each metric considered in this study. Results are illustrated for $\mathrm{NIW}=1,2,4$, and 8 (past) and $\mathrm{NOW}=1,2,$ and 4 (future). Metrics show consistent results (NRMSE and NAMMAE have to be minimized; R and AAM have to be maximized). Best predictions are achieved for the near future (short-term prediction, $\mathrm{NOW}=1$). Furthermore, the higher the NIW, the better is the prediction, even if a plateau is already achieved with $\mathrm{NIW}=4$. 

Figures \ref{fig:DMDder_5415M} and \ref{fig:DMDshift_5415M} extend the analysis to the ADMD formulation with time derivatives and time-shifted copies, respectively, focusing on the longer-term prediction results ($\mathrm{NOW}=4$), with respect to the $\mathrm{NIW}$, the number of derivatives ($\mathrm{NDE}$), and the number of time-shifted copies ($\mathrm{NTS}$). Since standard DMD has shown the consistency of the evaluation metrics, in the following, only NRMSE is used to describe the ADMD performance. Figure \ref{fig:DMDder_5415M} compares the standard DMD ($\mathrm{NDE}=0$) to its augmented counterpart with $\mathrm{NDE}=1, 2, 3$, and 4, providing the NRMSE box-plot. As for standard DMD, the higher the $\mathrm{NIW}$, the better are the results. In terms of time derivatives, the better results (on average) are achieved using $\mathrm{NDE}=2$, showing how further increasing $\mathrm{NDE}$ not necessarily improves ADMD performance. Figure \ref{fig:DMDshift_5415M} compares the standard DMD ($\mathrm{NTS}=0$) to its augmented counterparts with $\mathrm{NTS}=2, 4, 8$, and 16, providing the NRMSE box-plot. It may be noted that, since one encounter wave is composed by 32 time instant, the time-shifted copies under investigation correspond to $1/16$, $1/8$, $1/4$, and $1/2$ of one encounter wave. As for standard DMD and ADMD with derivatives, the higher the $\mathrm{NIW}$, the better are the results. In terms of time shifts, the better results (on average) seems to be function of $\mathrm{NIW}$, nevertheless (at least for the current analysis) the best trade-off is achieved using $\mathrm{NTS=4}$.  

\begin{figure*}[!t]
\centering
\includegraphics[width=1\textwidth]{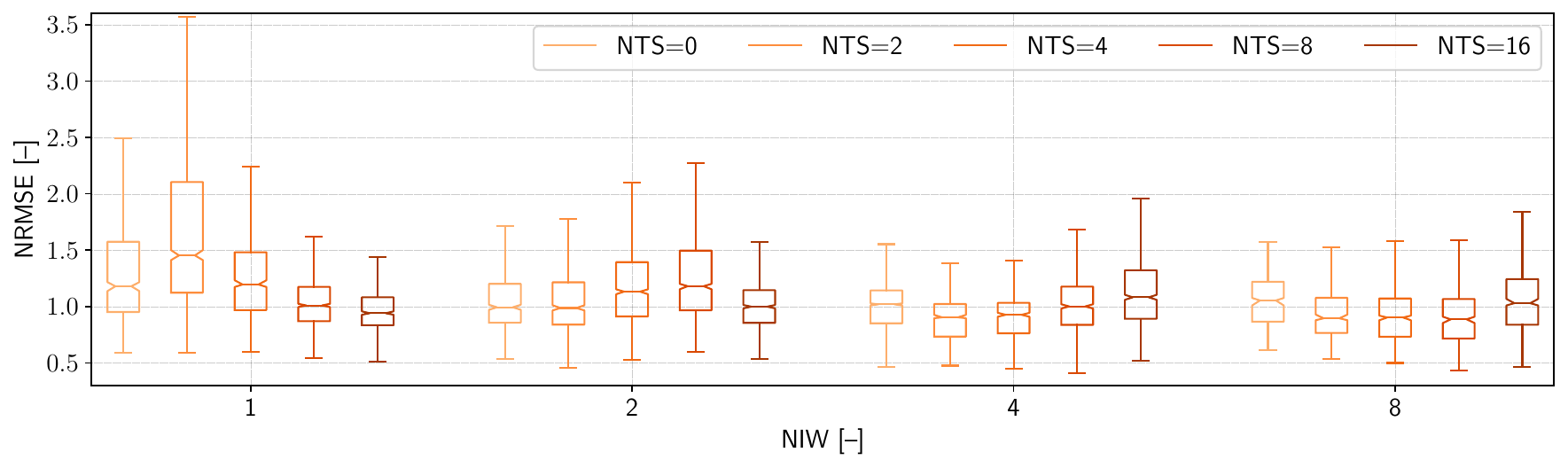}
\caption{Box plots of NRMSE values of ADMD with time-shifted copies for KCS test case.}\label{fig:DMDshift_KCS}
\end{figure*}
\begin{figure*}[!t]
\centering
\includegraphics[width=1\textwidth]{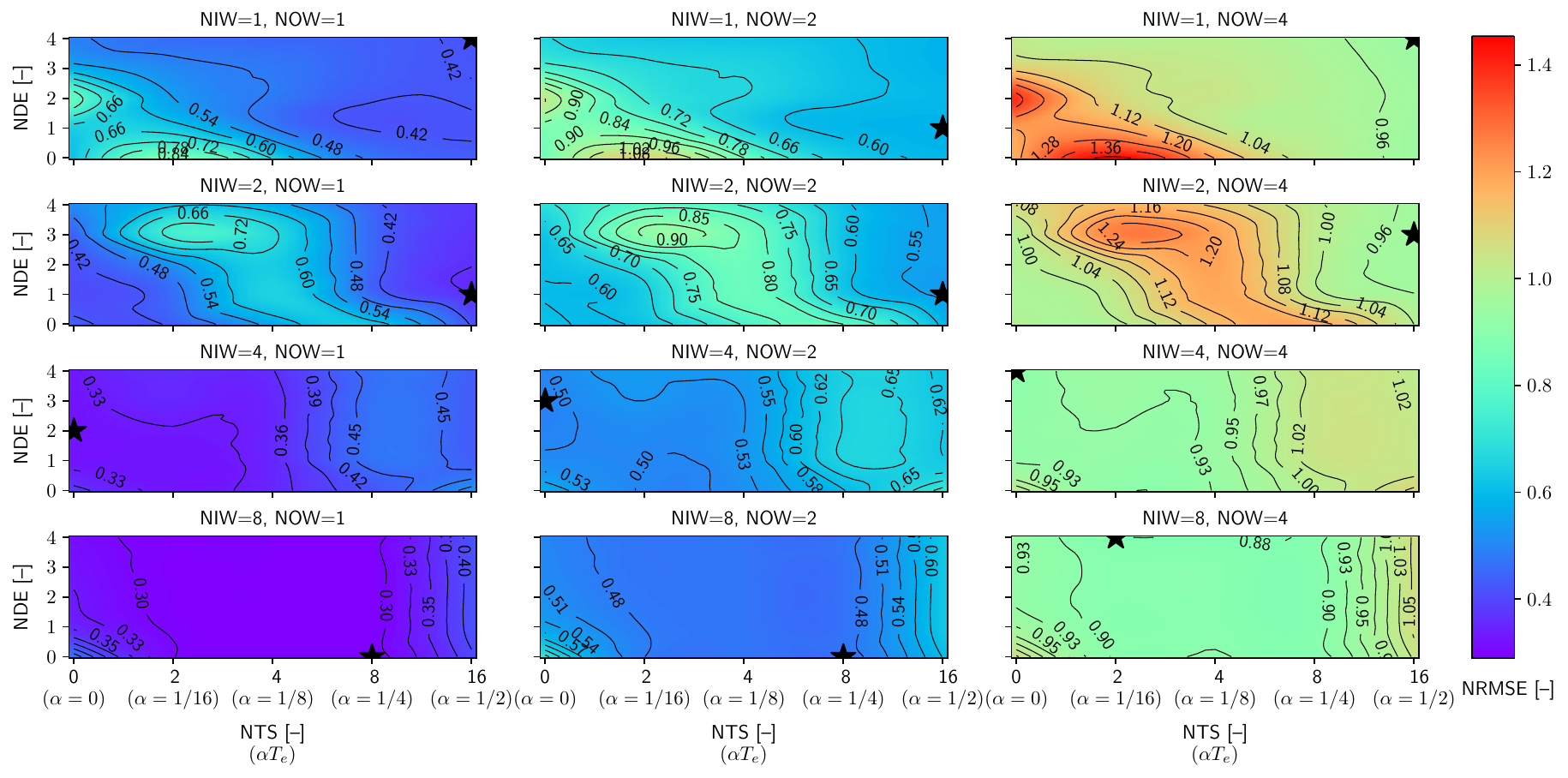}
\caption{Medians of NRMSE for the full factorial analysis of the DMD setups for KCS test case.}\label{fig:DMDcplot_KCS}
\end{figure*}
Figure \ref{fig:DMDcplot_5415M} gives the DMD/ADMD results of the full factorial combination of $\mathrm{NTS}$ and $\mathrm{NDE}$ ($x-$ and $y-$axis of each subplot, respectively) as a function of $\mathrm{NIW}$ (rows of Fig. \ref{fig:DMDcplot_5415M}) and $\mathrm{NOW}$ (columns of Fig. \ref{fig:DMDcplot_5415M}), where the contour plots represent the medians of NRMSE (blue is the best). On each subplot of Fig. \ref{fig:DMDcplot_5415M}, a black star is used to denote the best setup for the specific input/output combination of $\mathrm{NIW}$ and $\mathrm{NOW}$. Smallest NRMSE is achieved, as expected, for short-term prediction ($\mathrm{NOW}$) with the highest $\mathrm{NIW}$. On the contrary, larger errors are made for the prediction of higher $\mathrm{NOW}$ with a small $\mathrm{NIW}$ (small training set). A large $\mathrm{NTS}$ helps the future predictions when a low $\mathrm{NIW}$ is used, this because the use of time-shifted copies provides more information about the past. For what concerns the use of derivatives, it may be noted how intermediate values of $\mathrm{NDE}$ provide better results and their use always improves the performance of standard DMD ($\mathrm{NDE}=0$). Overall, for longer-term prediction ($\mathrm{NOW}=4$) the best results is achieved using $\mathrm{NIW}=8$, $\mathrm{NDE}=2$, and $\mathrm{NTS}=2$. This latter setup (ADMD) is finally compared to standard DMD in Fig. \ref{fig:DMDvsADMD_5415M}, where prediction of the system dynamics is shown for one random sample of the available data set. For the sake of clarity only the last 4 input waves are shown, depicted in black as observed (past) time histories, the predicted (future) time histories for DMD and ADMD are shown in blue and orange, respectively, while the true observed (future) time histories are presented with a dashed black line. All variables are shown in their standardized form and time values are normalized with the average encounter period ($T_e$). The NRMSE of the prediction is shown at the bottom. It may be noted how ADMD has evident better performance compared to standard DMD. In particular, ADMD provides better results for all variables except surge and heave. Even if the prediction is not perfect, the NRMSE shows how ADMD has always better future prediction than DMD. On average, variables are reasonably predicted up to $\mathrm{NOW}=2$ (two encounter periods). After, the prediction becomes less accurate, especially for DMD. 

\change{Finally, Figure \ref{fig:maxCorr_5415M} shows the relationship between the maximum value of the cross-correlation, $R_{xy}(\tau)$, between training and test sets (past and future time histories) and the corresponding DMD and ADMD prediction error (NRMSE). 
A negative correlation (as somehow expected) is found for heave, roll, pitch, yaw, and rudder angle, albeit rather mild. Conversely surge and sway show a mild positive correlation, which is to some extent counter intuitive. Although a general conclusion cannot be reached, based on the current results, on the impact of past/future correlation on the prediction capabilities, it may be noted how ADMD provides in general a slightly stronger negative correlation with smaller errors than DMD.} 

\begin{figure*}[!t]
\centering
\includegraphics[width=0.9\textwidth]{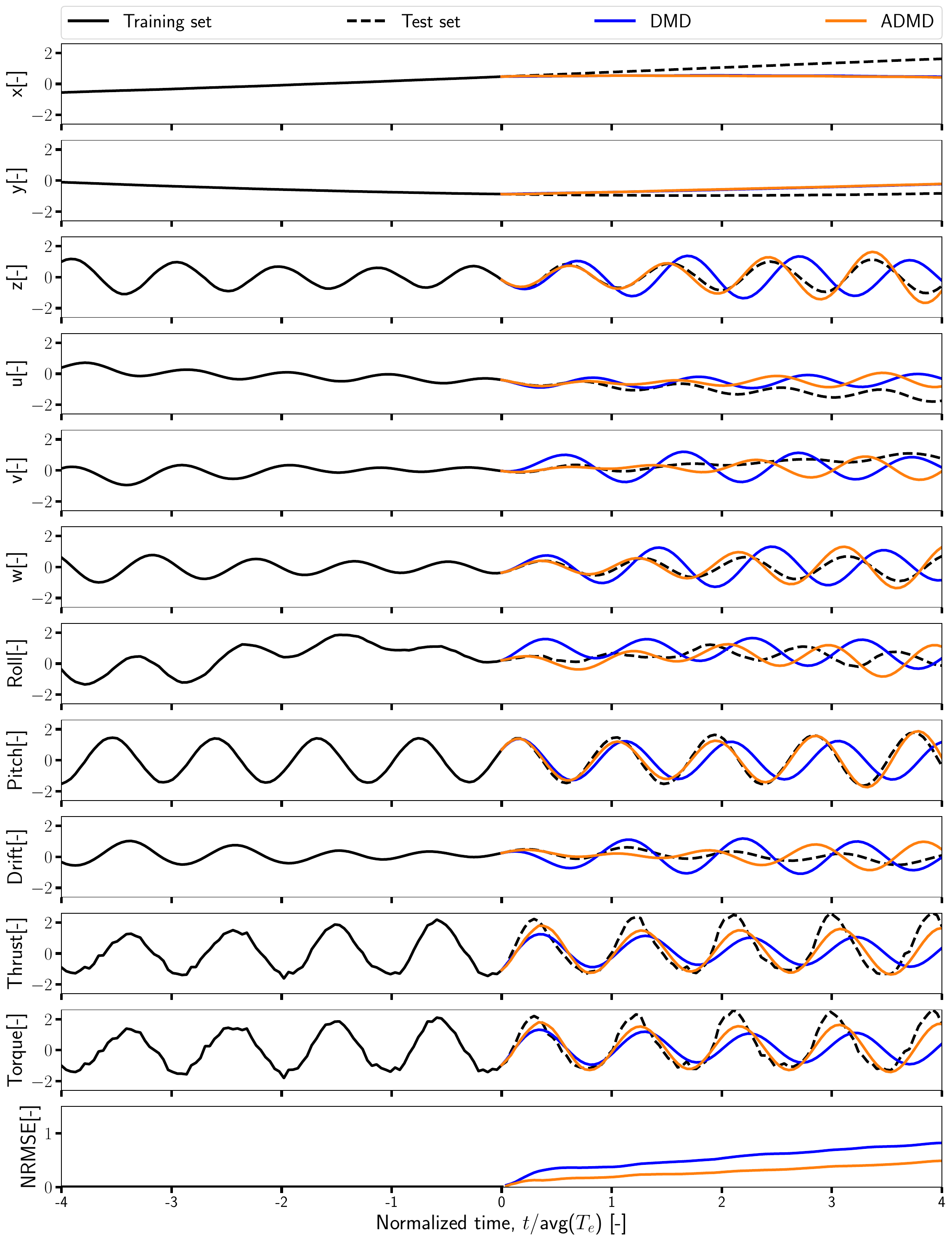}
\caption{Systems-state prediction comparison of DMD and best ADMD setup for a random time history of KCS data set.}\label{fig:DMDvsADMD_KCS}
\end{figure*}
\begin{figure*}[!t]
\centering
\includegraphics[width=1\textwidth]{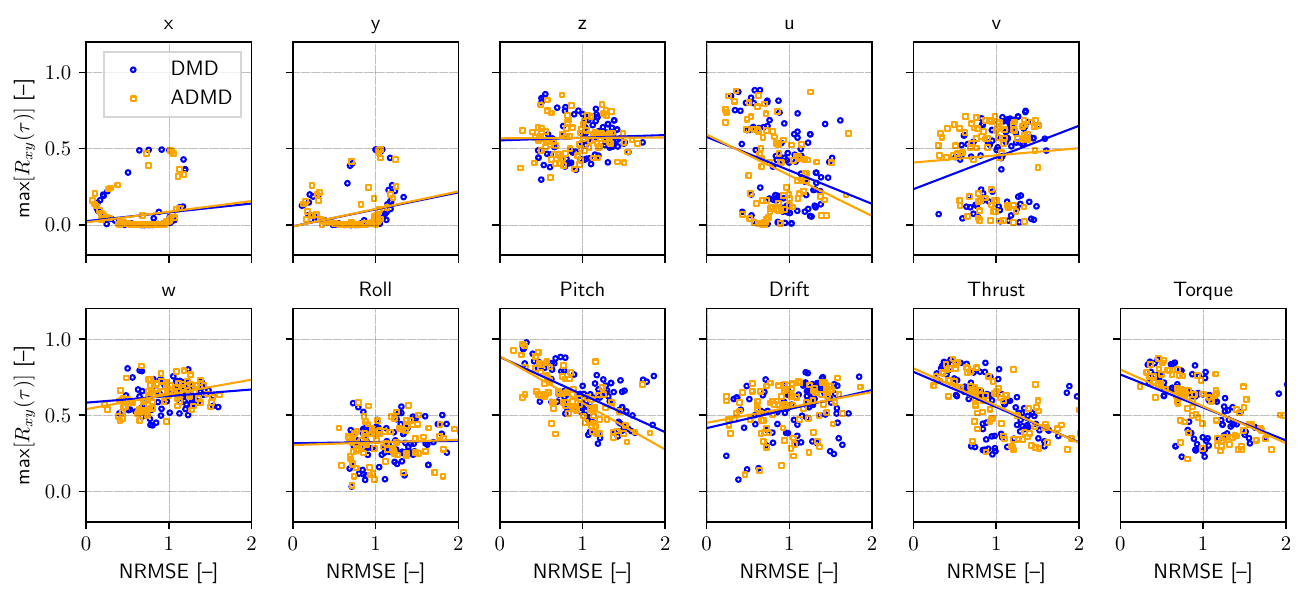}
\caption{Training/test sets maximum cross-correlation versus DMD and ADMD prediction error (NRMSE) for KCS case.}\label{fig:maxCorr_KCS}
\end{figure*}

\subsection{Forecasting of KCS Turning Circle in Regular Waves}
Figure \ref{fig:DMDstd_KCS} presents the prediction performance of the DMD standard formulation, showing the median values (of the statistic sample) for each evaluation metric. As for the 5415M case, the results are illustrated for $\mathrm{NIW}=1,2,4$, and 8 (past) and $\mathrm{NOW}=1,2$, and 4 (future). Metrics show consistent results, expect for NAMMAE with $\mathrm{NOW}=4$. Differently from the 5415M case, best predictions are achieved with an intermediate $\mathrm{NIW}$ for both short- and longer-term prediction ($\mathrm{NIW}=4$ for $\mathrm{NOW}=1$ and 2 and $\mathrm{NIW}=2$ for $\mathrm{NOW}=4$, respectively). This can be related to the specific dynamics that is observed, since during turning circle in waves, the boundary conditions (incident wave direction) change with the maneuvers trajectory, so during future time instant, far from the past observation, the observed dynamic can not be longer related with the past.

Figures \ref{fig:DMDder_KCS} and \ref{fig:DMDshift_KCS}, as for the 5415M case, extend the analysis to the ADMD formulation with time derivatives and time-shifted copies, respectively, focusing on the longer-term prediction results ($\mathrm{NOW}=4$). Also in this case, only NRMSE is used to describe the ADMD performance. Figure \ref{fig:DMDder_KCS} compares the standard DMD ($\mathrm{NDE}=0$) to its augmented counterparts with $\mathrm{NDE}=1, 2, 3$, and 4, providing the NRMSE box-plot. Contrary to standard DMD, the higher the $\mathrm{NIW}$, the better are the results. In terms of time derivatives, the better results (on average) are achieved using $\mathrm{NDE}=2$, where a plateau is reached for higher $\mathrm{NIW}$. Figure \ref{fig:DMDshift_KCS} compares the standard DMD ($\mathrm{NTS}=0$) to its augmented counterpart with $\mathrm{NTS}=2, 4, 8$, and 16, providing the NRMSE box-plot. As for standard ADMD with derivatives, the higher the $\mathrm{NIW}$, the better are the results. In terms of time shifts, the results are consistent with those obtained for the 5415M case: the better results (on average) seems to be function of $\mathrm{NIW}$, nevertheless (at least for the current analysis) the best trade-off is achieved using $\mathrm{NTS=4}$. 

Figure \ref{fig:DMDcplot_KCS}, with the same layout described for the 5415M case, gives the DMD/ADMD results of the full-factorial combination of $\mathrm{NTS}$ and $\mathrm{NDE}$ as a function of $\mathrm{NIW}$ and $\mathrm{NOW}$, where the contour plots represent the medians of NRMSE. Results obtained are quite consistent with the 5415M analysis, providing similar trends, even if the errors are generally higher. Smallest NRMSE is achieved, for limited $\mathrm{NOW}$ with the highest $\mathrm{NIW}$. On the contrary, larger errors are made for the prediction of higher $\mathrm{NOW}$ with a small $\mathrm{NIW}$. A large $\mathrm{NTS}$ helps the future predictions when a low $\mathrm{NIW}$ is used. For what concerns the use of derivatives, no specific trends can be highlighted. Overall, for longer-term prediction ($\mathrm{NOW}=4$) the best results is achieved using $\mathrm{NIW}=8$, $\mathrm{NDE}=4$, and $\mathrm{NTS}=2$. This latter setup (ADMD) is finally compared to standard DMD in Fig. \ref{fig:DMDvsADMD_KCS}, where prediction of the system dynamics is shown for one random sample of the available data set, as for the 5415M case. It may be noted how ADMD has better performance compared to standard DMD. In particular, ADMD provides visibly better results for $z$, \change{$w$, pitch,} thrust, and torque. Also in this case prediction degrades with time, nevertheless, the NRMSE shows how ADMD has always better future prediction than DMD.

\change{Finally, Figure \ref{fig:maxCorr_KCS} shows the relationship between the maximum value of the cross-correlation, $R_{xy}(\tau)$, between training and test sets (past and future time histories) and the corresponding DMD and ADMD prediction error (NRMSE). In this case, a negative correlation (expected) is found for $u$, pitch, thrust, and torque, albeit again rather mild. Conversely, $x$, $y$, $w$, and drift show a mild positive correlation, which is again counter intuitive; $z$ and roll exhibit a nearly zero correlation. Also in this case, although a general conclusion cannot be reached on the impact of past/future correlation on the prediction capabilities, it is shown how the ADMD provides in general a slightly stronger negative correlation with smaller errors than DMD.}

\subsection{Modal Analysis of 5415M and KCS}
\begin{figure*}[!t]
\centering
\includegraphics[width=1\textwidth]{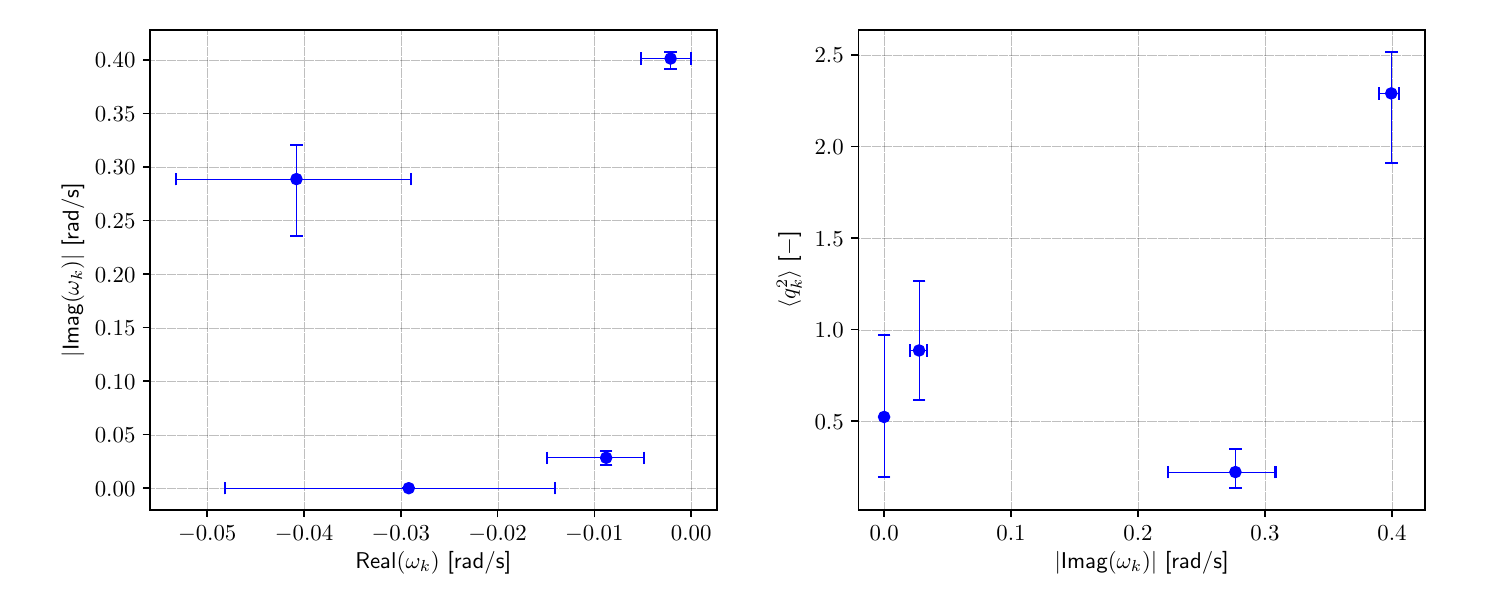}
\caption{DMD complex modal frequencies (left) and modal participation (right) for the 5415M test case.}\label{fig:DMDmodes_5415M}
\end{figure*}
\begin{figure}[!t]
\centering
\includegraphics[width=0.5\columnwidth]{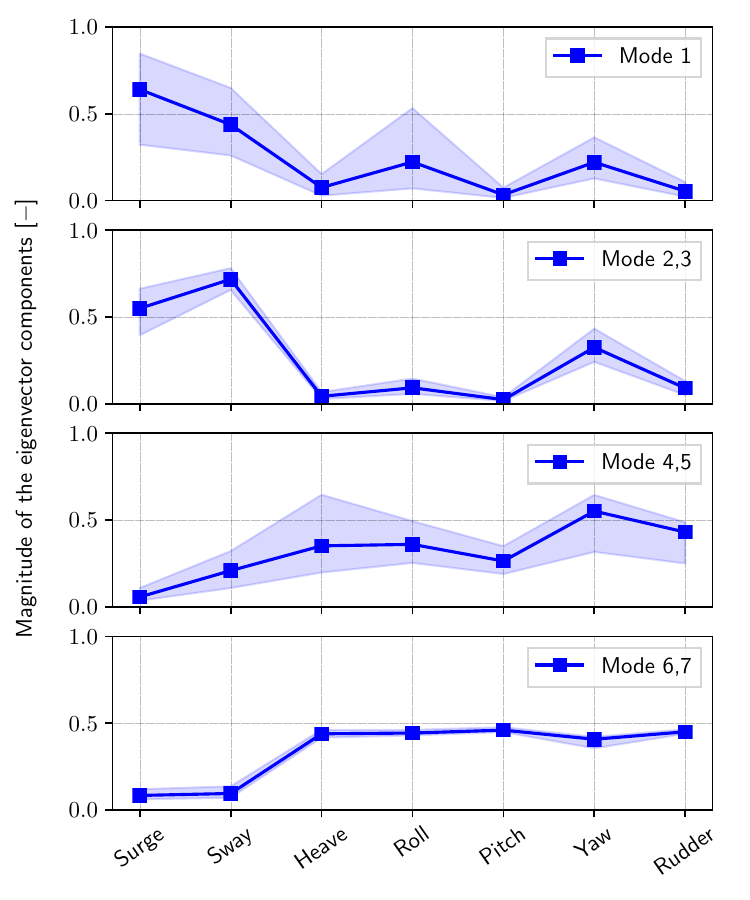}
\caption{The eigenvector components of all modes sorted by frequency (in ascending order) for the 5415M test case.}\label{fig:DMD12modes_5415M}
\end{figure}
\change{
A modal analysis is carried out for both the 5415M and the KCS following a statistical approach similar to that used for the assessment of the DMD forecasting capabilities. Specifically, 1001 time windows, each including 8 encounter waves on average, are randomly selected from the available time histories. The DMD is performed for each time window, where eigenvalues and  eigenvectors are ordered by the absolute value of the imaginary part of the associated complex frequency (the frequency in a strict sense), in ascending order. For each mode, the modal participation (or energy) is assessed by the mean square of the modal coordinate, $\langle q^2 \rangle$, evaluated by projection of the data onto the mode. Complex-conjugate modes are coupled together and assessed as one, as far as the modal participation is concerned.}

\change{
Figure \ref{fig:DMDmodes_5415M} (left) shows the complex frequencies provided by the DMD. Median and uncertainty bands provided by the inter-quartile range are depicted. Figure \ref{fig:DMDmodes_5415M} (right) shows the modal participation as a function of the frequency. The corresponding modes are shown in Fig. \ref{fig:DMD12modes_5415M}. The median of the magnitude of the complex-conjugate vector component is depicted along with the uncertainty band provided by the inter-quartile range. It may be noted how the most energetic dynamics is found close to 0.40 rad/s, with a very narrow uncertainty band along the frequency axis. This corresponds to the response of the ship to the peak frequency of the encounter spectrum, which equals 0.39 rad/s and is close to the resonance condition for the roll motion \cite{serani2021-OE}. The corresponding mode is shown at the bottom of Fig. \ref{fig:DMD12modes_5415M} revealing a significant participation and coupling of heave, roll, pitch, sway motions and rudder angle. Planar motion variables (surge, sway, yaw) coupled with rudder angle are associated to the second most energetic mode, which exhibits a significantly slower dynamics.
}

\begin{figure*}[!t]
\centering
\includegraphics[width=1\textwidth]{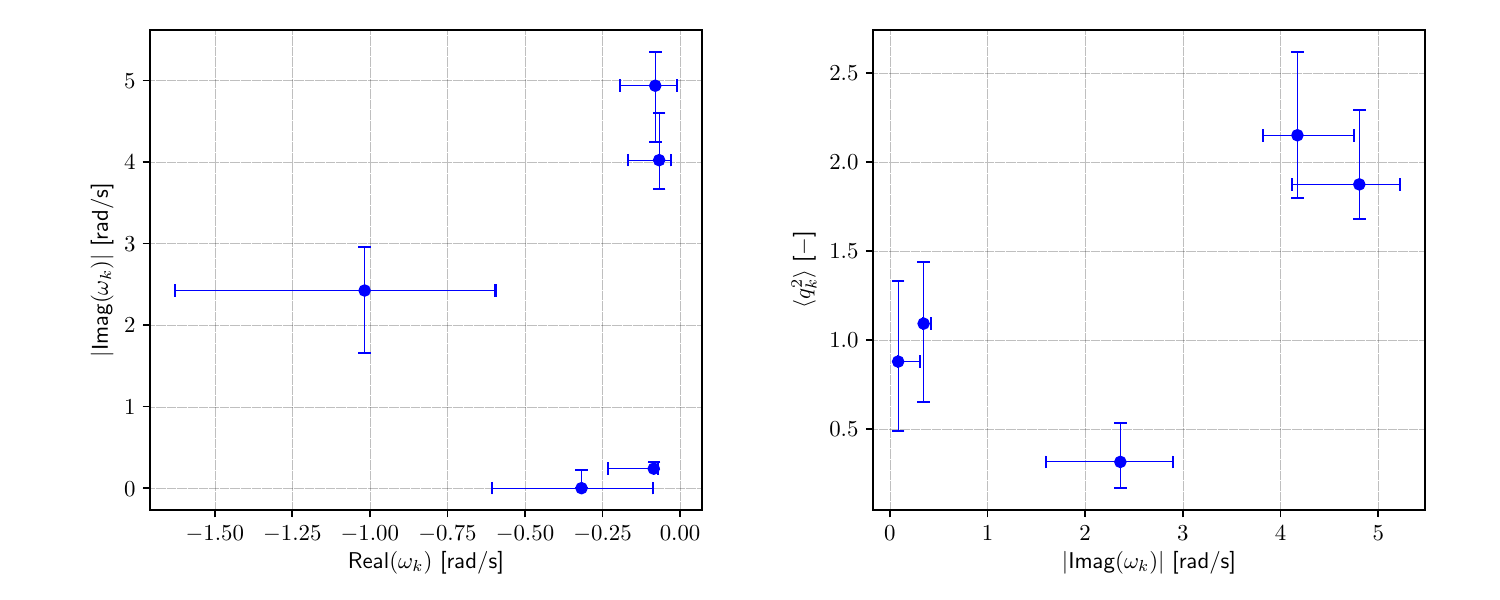}
\caption{DMD complex modal frequencies (left) and modal participation (right) of all modes for the KCS test case.}\label{fig:DMDmodes_KCS}
\end{figure*}

\begin{figure}[!t]
\centering
\includegraphics[width=0.5\columnwidth]{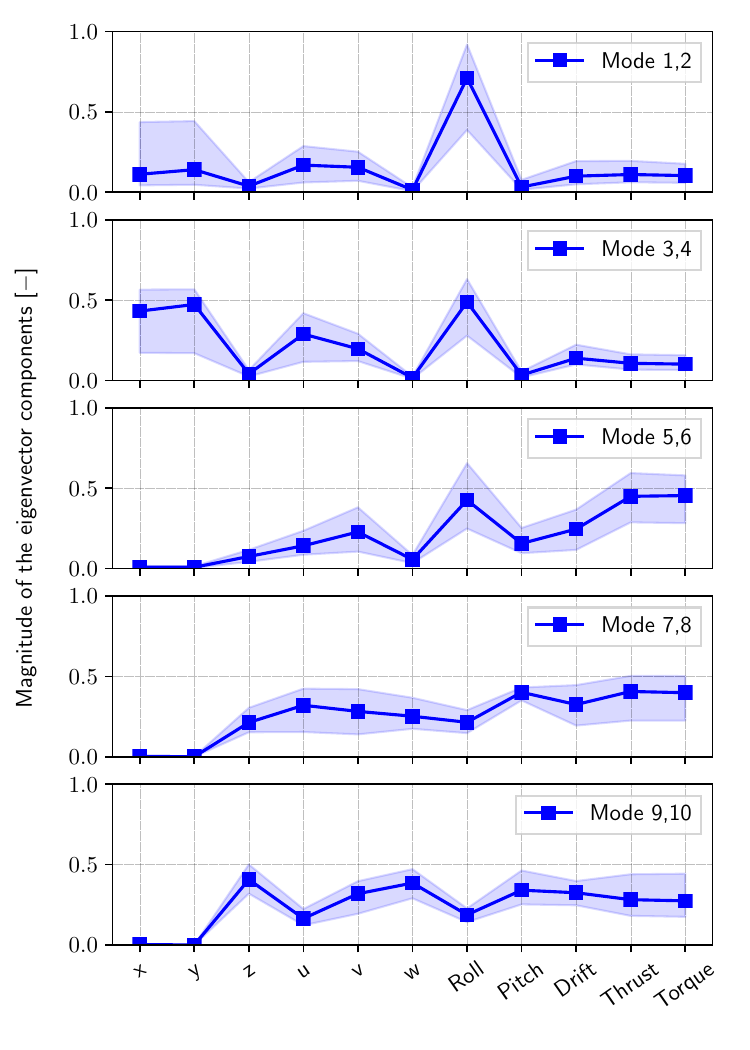}
\caption{The eigenvector components of all modes sorted by frequency (in ascending order) for the KCS test case.}\label{fig:DMD12modes_KCS}
\end{figure}

\change{Similar results are shown for the KCS in Figs. \ref{fig:DMDmodes_KCS} and \ref{fig:DMD12modes_KCS}. The most energetic modes fall in the frequency range associated to the encounter frequencies for head and beam waves, showing participation of heave ($z$), pitch, surge ($u$), along with significant coupling with all other variables but trajectory variables $x$ and $y$. The latter participate in a significantly slower mode associated to accomplishing the turning circle. Finally, it may be noted how the roll motion participates both in the slow turning-circle dynamics and in a mode that falls in the 1.5--3.0 rad/s frequency range, which corresponds to the ship response to following waves. This latter mode is also significantly participated by thrust and torque.
}

%
%

%



\section{Conclusions and Future Work}\label{sec:conclusions}
A statistical analysis on the use of dynamic mode decomposition and its augmented variant to forecast  trajectories, motions, and forces of ships operating in waves was presented and discussed. The DMD provides a linear representation of the system dynamics, allowing for (a) forecasting the system's state in the near future and (b) extracting knowledge on the system dynamics. Results were shown for course keeping data of the self-propelled 5415M in irregular waves and turning-circle data of the self-propelled KCS in regular waves. Time histories were provided by CFD and EFD for 5415M and KCS, respectively.

The statistical assessment was based on four evaluation metrics varying the number of input (training/past states) and output (test/future states prediction) waves, as well as exploiting the state augmentation, by adding derivatives and shifted copies of the time histories. A full-factorial combination of these parameters has been investigated.  

Overall state augmentation strategy (ADMD) combining both time derivatives and time-shifted copies showed better prediction capabilities for both test cases, compared to standard DMD. In general, the longer is the time horizon to be predicted, the longer is the past-observation window needed to feed the method and get an accurate forecast. The use of time derivatives showed that using up to the second derivatives improves the quality of the forecast. Similarly, the use of time-delayed copies of time histories provided the best results using a delay up to 1/8 of the encounter wave period. Overall, for longer term prediction ($\mathrm{NOW}=4$) the combination of $\mathrm{NDE}=2$ and $\mathrm{NTS}=2$ and 4 (for 5415M and KCS, respectively) provided the most promising results. The analysis is very efficient and suitable for real-time predictions. 

In conclusion, the method is able to identify the most important modes in the observed dynamics and forecast, with reasonable accuracy, the system's state up to two wave encounter periods. After this time horizon, the prediction is no longer accurate and the methodology needs further improvements.

Methodological advancements that will be subject of future work include %
DMD with time-delay embedding \cite{kamb2020time} and multi-resolution DMD \cite{kutz2016multiresolution}. Finally, the combination of DMD with artificial NN approaches \cite{dagostino2022-OEME} is expected to overcome some of the limitations of the DMD (i.e., its linearity) providing more flexible architectures to address highly-nonlinear system dynamics. For this reason, on going research activities are focusing on the development of hybrid architectures based on DMD and RNNs \cite{diez2022snh}. The assessment of DMD prediction performance conditional to different operating and environmental conditions is also of interest and will be assessed within the activity of the NATO AVT-348.
Finally, the current method will be extended to system-identification problems, via DMD with control \cite{proctor2016dynamic}. 

\section*{Acknowledgments}
CNR-INM is grateful to Drs. Elena McCarthy and Woei-Min Lin of the Office of Naval Research for for their support through the Naval International Cooperative Opportunities in Science and Technology Program, grant N62909-21-1-2042.
Dr. Andrea Serani is also grateful to the National Research Council of Italy, for its support through the Short-Term Mobility Program 2018. 
The research is conducted in collaboration with NATO STO Research Task Groups AVT-280 ``Evaluation of Prediction Methods for Ship Performance in Heavy Weather,''  AVT-348 ``Assessment of Experiments and Prediction Methods for Naval Ships Maneuvering in Waves,'' and AVT-351 ``Enhanced Computational Performance and Stability \& Control Prediction for NATO Military Vehicles''.

\bibliographystyle{abbrv}  
\bibliography{biblio}  

\end{document}